\documentclass[12pt]{article}
\usepackage[american]{babel}
\usepackage{amsmath}
\usepackage{latexsym}
\usepackage{amssymb}
\usepackage{theorem}
\usepackage{diagrams}
\usepackage{mathrsfs}
\usepackage{calligra}
\DeclareMathAlphabet{\mathcalligra}{T1}{calligra}{m}{n}
\DeclareFontShape{T1}{calligra}{m}{n}{<->s*[1.2]callig15}{}
\theoremstyle{change}
\newtheorem{Thm}{Theorem}[section]
\newtheorem{Cor}[Thm]{Corollary}
\newtheorem{Prop}[Thm]{Proposition}
\newtheorem{Lem}[Thm]{Lemma}
{\theorembodyfont{\rmfamily}
\newtheorem{Num}[Thm]{}

}
\newcommand{\Ind}{\mathrm{Ind}}
\newcommand{\ind}{\mathrm{ind}}
\renewcommand{\phi}{\varphi}

\newcommand{\bra}[1]{\langle#1\rangle}
\newcommand{\proof}{\par\medskip\rm\emph{Proof. }}

\newcommand{\qed}{\ \hglue 0pt plus 1filll $\Box$}

\newcommand{\mapstoo}{\longmapsto}

\newcommand{\RR}{\mathbb{R}}
\newcommand{\ZZ}{\mathbb{Z}}

\renewcommand{\SS}{\mathbb{S}}

\newcommand{\SKIP}[1]{}

\newcommand{\eps}{\varepsilon}
\renewcommand{\emptyset}{\varnothing}

\newcommand{\supp}{\mathrm{supp}}

\renewcommand{\setminus}{-}
\begin{document}

\title{{\bf On the local structure and the homology\\ of CAT$(\kappa)$ spaces and euclidean buildings}}
\author{Linus Kramer\thanks{This work was supported by SFB 878.}\\ \\
\textit{Dedicated to Helmut R. Salzmann}}
\maketitle
\begin{abstract}
We prove that every open subset of a euclidean building is
a finite dimensional absolute neighborhood retract. This implies in particular
that such a set has the homotopy type of a finite dimensional simplicial complex. We also include
a proof for the rigidity of homeomorphisms of euclidean buildings.
A key step in our approach
to this result is the following: the space of directions $\Sigma_oX$ of
a CAT$(\kappa)$ space $X$ is homotopy equivalent to a small
punctured disk $B_\eps(X,o)\setminus o$. The second ingredient is the local
homology sheaf of $X$. Along the way, we prove some results
about the local structure of CAT$(\kappa)$-spaces which may be of
independent interest. 
\end{abstract}

A CAT$(\kappa)$ space $X$ is a metric space
where geodesic triangles of perimeter less than $2D_\kappa$
exist and are not thicker than in the comparison space $M_\kappa$ of constant sectional
curvature $\kappa$ and diameter $D_\kappa$ (where $D_\kappa=\pi/\sqrt\kappa$ for $\kappa>0$
and $D_\kappa=\infty$ for $\kappa\leq 0$). Bridson-Haefliger \cite{BH} and 
Burago-Burago-Ivanov \cite{BBI} give a thorough introductions to these spaces.
More generally, a metric space is said to have curvature bounded above if every point
has a neighborhood which is CAT$(\kappa)$, for some
$\kappa$.
We prove the following results.

\medskip\noindent
\textbf{Theorem A}
{\em Let $X$ be a metric space with curvature bounded above.
For every point $o\in X$, there is a number $\eps>0$ such that the punctured ball
$B_\eps(o)\setminus o$ is homotopy equivalent to the space of
directions $\Sigma_oX$.}

\medskip\noindent
The proof uses at some point that spaces with curvature bounded above are
ANRs. See Theorem~\ref{CATImpliesANR} and the remarks preceding it.

\medskip\noindent
\textbf{Theorem B}
{\em Let $X$ be a euclidean building with $n$-dimensional apartments.
Then $X$ is an $n$-dimensional absolute retract. 
Every nonempty open subset of $X$ is an $n$-dimensional absolute neighborhood retract.}

\medskip\noindent
As an application, we give a proof of the following result.

\medskip\noindent
\textbf{Theorem C}
{\em Let $f:X\rTo X'$ be a homeomorphism between metrically irreducible euclidean
buildings of dimension at least $2$.
Then $f$ maps apartments in the maximal atlas of $X$ to apartments in the maximal atlas of $X'$.
Possibly after rescaling $X'$, there exists a unique isometry $\tilde f:X\rTo X'$
with $f(A)=\tilde f(A)$ for all apartments $A\subseteq X$. If $X$ has more than
one thick point, then $f$ has finite distance from $\tilde f$.}

\medskip
\medskip
After the paper was completed, I learned that Theorem~A had also been obtained
by Lytchak and Nagano some time ago; see \cite{LytNag}. Bruce Kleiner informed
me that he had also obtained this result (unpublished).
Also, Lang and Schlichenmaier had computed the Nagata
dimension of euclidean buildings in \cite{LaS}. Their result implies Theorem~B.
The result that CAT$(\kappa)$ spaces
are ANRs has been obtained by various authors (in different degrees of generality).
Theorem~C was proved by Kleiner
and Leeb in \cite{KL} for the case of complete euclidean buildings.
More references and details are given below.
I thank Mladen Bestvina, 
Martin Bridson, Jerzy Dydak, Misha Kapovich, Bruce Kleiner,
Urs Lang, Alexander Lytchak,
Sibe Marde{\v{s}}i{\'c} and Stephan Stolz for helpful remarks.

\medskip
I learned geometry and topology (and some dimension theory) first from Reiner Salzmann.
His concise (and often demanding) lectures in T\"ubingen were among
the finest math classes that I ever took.
It is a pleasure to dedicate this article to him on the occasion of
his eightieth birthday.

\medskip\noindent\textbf{Outline of the paper.}\\
In Section~1 we collect some general facts about spaces with upper curvature bounds.
We introduce in Section~2 the space of directions and the 'blow-up' at a given point.
Section 3 contains a proof of Theorem~A. Our approach is essentially homotopy theoretic.
Some simple sheaf theoretic notions are introduced in Section~4. In Section~5, these
are applied to the homology of spherical buildings. Section~6 combines the previous
results into a proof
of Theorem~C. In Section~7, we prove Theorem~B. The Appendix serves as a reference
for some topological results on ANRs and the covering dimension.

The definition of a spherical building that we use is 'simplicial'
as in Tits' Lecture Note \cite[3.1]{TitsLNM}
or \cite[1.1]{TitsComo},
except that we allow weak (non-thick) spherical buildings; see
\ref{SphericalBuildings} below. We define euclidean buildings essentially as Tits does in
\cite[1.4]{TitsComo}, see also \cite{KL} and in particular \cite{Par}. Euclidean buildings 
are CAT$(0)$ spaces which can be highly
singular (they need not be simplicial and they may branch everywhere);
see \ref{EBAxioms}.

\section{Spaces with curvature bounded above}

We first introduce some notation.
Suppose that $(X,d)$ is a  metric space, that $o\in X$ and that $r>0$.
We put
\[
B_r(o)=\{x\in X\mid d(x,o)<r\}\text{ and }
\bar B_r(o)=\{x\in X\mid d(x,o)\leq r\}.
\]
We say that $X$ has \emph{radius at most $r$} if $X\subseteq \bar B_r(o)$ 
for some $o\in X$, and that $X$ has \emph{diameter at most $r$} if
$X\subseteq \bar B_r(o)$ for every $o\in X$.
Two maps $f,g:Z\pile{\rTo\\ \rTo} X$ are \emph{$r$-near} if
$d(f(z),g(z))\leq r$ holds for all $z\in Z$.
A \emph{geodesic} in $X$ is an isometry $J\rTo X$, 
where $J\subseteq \RR$ is a closed interval. We call
the image of $J$ a \emph{geodesic segment} in $X$.

We now collect some facts about spaces with curvature bounded above. 
The books by Bridson-Haefliger \cite{BH} and Burago-Burago-Ivanov \cite{BBI}
are excellent introductions to this area of geometry.
\begin{Num}\textbf{The model spaces $M_\kappa$.}\label{ModelSpaces}\\ 
The \emph{model spaces} $M_{-1}$, $M_0$ and $M_1$ are the
$2$-dimensional
\footnote{Sometimes it is convenient to have higher dimensional model
spaces. For our purposes the surfaces suffice.} 
hyperbolic space $\mathbb H^2$, euclidean space $\mathbb E^2$ and the
$2$-sphere $\SS^2$, respectively. These spaces are endowed with the standard
hyperbolic, euclidean and spherical metrics $d_H,d_E,d_S$, respectively;
see \cite[I.2]{BH}.
For general $\kappa<0$ the model space $M_\kappa$ is defined
as $(M_{-1},\sqrt{-\kappa}d_H)$, and for $\kappa>0$ as
$(M_1,\sqrt\kappa d_S)$. These are precisely the complete simply connected
surfaces of constant sectional curvature $\kappa$ and diameter $D_\kappa$,
where $D_\kappa=\pi/\sqrt\kappa$ for $\kappa>0$
and $D_\kappa=\infty$ for $\kappa\leq 0$.
\end{Num}
\begin{Num}\textbf{CAT$(\kappa)$ spaces.}\label{CAT(k)}\\
Let $\kappa$ be a real number and let $M_\kappa$ and $D_\kappa$ be as in
\ref{ModelSpaces}. A metric space is called a \emph{CAT$(\kappa)$ space}
if any two points at distance less than $D_\kappa$ can be joined by
a geodesic, and if geodesic triangles of perimeter less than
$2D_\kappa$ are not thicker than the comparison triangles in $M_\kappa$;
see \cite[II.1.1]{BH} for more details. This condition on triangles
is called the \emph{CAT$(\kappa)$ inequality}. 

The metric completion of
a CAT$(\kappa)$ space is again CAT$(\kappa)$; see \cite[II.3.11]{BH}. If $X$ is CAT$(\kappa)$,
then $X$ is also CAT$(\kappa')$, for all $\kappa'\leq\kappa$; see \cite[II.1.12]{BH}.
A metric space has \emph{curvature bounded above} if every point has a neighborhood
which is CAT$(\kappa)$, for some $\kappa$ (which may depend on the point).

The CAT$(\kappa)$ inequality implies that points $u,v$ at distance 
$d(u,v)<D_\kappa$ can be joined by a unique geodesic. The corresponding
geodesic segment is denoted $[u,v]$. It varies continuously with $u$ and $v$;
see \cite[II.1.4]{BH}.
If $Z$ is a topological space and if $f,g:Z\pile{\rTo\\ \rTo} X$ are continuous
and $r$-near, for some $r<D_\kappa$, then
$f$ and $g$ are homotopic via the \emph{geodesic homotopy} that
moves $f(z)$ along $[f(z),g(z)]$ to $g(z)$.
\end{Num}
\begin{Num}\textbf{Convex sets and convex hulls}\label{Convex}\\
A subspace $C$ of a CAT$(\kappa)$ space of diameter less than $D_\kappa$
is \emph{convex} if $[u,v]\subseteq C$ holds for
all $u,v\in C$. For $o\in X$, the balls $B_r(o)$ and $\bar B_r(o)$ are convex
if $r< D_\kappa/2$; see \cite[II.1.4]{BH}. If $A\subseteq X$ has radius 
$r<D_\kappa/2$, then the \emph{convex hull} of $A$ is defined
to be the intersection of all convex sets containing $A$.
The convex hull contains $A$, is convex, CAT$(\kappa)$, and has radius at most $r$.
A convex set of radius less than $D_\kappa/2$ is contractible
by a geodesic homotopy; see \cite[II.1.5]{BH}.
\end{Num}
\begin{Num}\textbf{The Alexandrov angle}\label{Angles}\\
Suppose that $o$ is a point in the CAT$(\kappa)$ space $X$ and that
$c:[0,a]\rTo X$ and $c':[0,a']\rTo X$ are two geodesics 
starting at $o$. The \emph{Alexandrov angle} between $c$ and $c'$ is
\[
\angle(c,c')=2\lim_{s\to0}\arcsin(d(c(s),c'(s))/2s);
\]
see \cite[II.3.1]{BH}
\footnote{The Alexandrov angle in $M_\kappa$ coincides with
the usual Riemannian angle.}. 
For the endpoints $u=c(a)$ and $v=c'(a')$ we put
\[
\angle_o(u,v)=\angle(c,c').
\]
The angle $(u,v)\mapstoo\angle_o(u,v)$ is a continuous
pseudometric on $B_{D_\kappa}(o)\setminus o$; see \cite[I.1.14,II.3.3]{BH}.
The Alexandrov angles in a triangle in $X$ are not greater than the angles
in the comparison triangle in $M_\kappa$; see \cite[II.1.7(4)]{BH}.
\end{Num}
\begin{Num}\textbf{The Law of Sines in $M_\kappa$}\label{LawOfSines}\\
Let $a,b,c$ be a triangle of diameter less than $2D_\kappa$ in $M_\kappa$.
Let $\alpha,\beta,\gamma$ denote the angles at $a,b,c$. We will use
the following two formulas.

\medskip
\noindent\textbf{The euclidean case $\kappa=0$.}
If $\kappa=0$, then we have the well-known euclidean \emph{Law of Sines:}
\[
\frac{\sin(\alpha)}{d(b,c)}=
\frac{\sin(\beta)}{d(c,a)}=
\frac{\sin(\gamma)}{d(a,b)}.
\]

\medskip
\noindent\textbf{The case $\kappa>0$.}
For $\kappa>0$, the model space $M_\kappa$ is the unit sphere $\SS^2$ in
$\mathbb E^3$ with metric $d(u,v)=\frac1{\sqrt\kappa}\arccos(\bra{u,v})$.
The spherical \emph{Law of Sines} in $M_\kappa$ is then
\[
\frac{\sin(\alpha)}{\sin(\sqrt\kappa d(b,c))}=
\frac{\sin(\beta)}{\sin(\sqrt\kappa d(c,a))}=
\frac{\sin(\gamma)}{\sin(\sqrt\kappa d(a,b))}.
\]
\proof
To see this, let $u,v$ be unit tangent vectors at $a$ pointing towards
$b$ and $c$, respectively. Then $u\times v=a\sin(\alpha)$
\footnote{We may assume that $a,u,v$ is right-handed.}, whence
\[
\det(a,u,v)=\bra{a,u\times v}=\sin(\alpha).
\]
We have $b=a\cos(\sqrt\kappa d(a,b))+u\sin(\sqrt\kappa d(a,b))$ and
$c=a\cos(\sqrt\kappa d(a,c))+v\sin(\sqrt\kappa d(a,c))$.
Thus
\begin{align*}
\det(a,b,c)&=\bra{a,b\times c}\\
&=\bra{a,u\times v} \sin(\sqrt\kappa d(a,b))\sin(\sqrt\kappa d(a,c)) \\
&=\sin(\alpha)\sin(\sqrt\kappa d(a,b))\sin(\sqrt\kappa d(a,c)).
\end{align*}
Permuting $a,b,c$ cyclically we obtain three times the same value and
the formula above follows.
\qed
\end{Num}

\section{The space of directions and the blow-up at a point}
\label{Section2}
In this section we introduce the space of direction $\Sigma_oX$ at
a point (the 'unit tangent space' at $o$). We also show that
by a change of the metric,
$\Sigma_oX$ can be glued into the puncture $o$ of $X\setminus o$.
We call the resulting space the blow-up of $X$ at $o$.

The following terminology will be fixed throughout this section.
We assume that $X$ is a CAT$(\kappa)$ space with $o\in X$ such that
\[
X\subseteq \bar B_{D_\kappa/2}(o).
\]
We put $Y=X\setminus o$. For $x\in X$ and $s\in [0,1]$, we denote by
$sx$ the unique point in the geodesic segment $[o,x]$ with $d(sx,o)=sd(x,o)$.
\begin{Num}\textbf{The space of directions.}\label{SigmaX}\\
We noted in \ref{Angles} that the Alexandrov angle induces
a pseudometric $(u,v)\mapstoo\angle_o(u,v)$ on $Y$. 
The metric completion of $Y$ with respect to $\angle_o$ is the
\emph{space of directions} $\Sigma_oX$ (where
points $x,y$ with $\angle_o(x,y)=0$ are identified);
see \cite[II.3.18 and 19]{BH}. The space
$(\Sigma_oX,\angle_o)$ is always CAT$(1)$; see \cite[II.3.19]{BH}.
We denote the canonical map $Y\rTo\Sigma_0X$ by $\rho$.

In order to make this map Lipschitz, we change the metric on $Y$.
For $x,y\in Y$ we put
\[
d_o(x,y)=\sqrt{d(x,y)^2+\angle_o(x,y)^2}.
\]
The identity map $(Y,d_o)\rTo(Y,d)$ is obviously
$1$-Lipschitz, and so is $\rho:(Y,d_o)\rTo(\Sigma_oX,\angle_o)$.
\end{Num}
\begin{Lem}
The identity $(Y,d)\rTo(Y,d_o)$ is locally uniformly continuous.

\proof
Let $p\in Y$.
Suppose first that $\kappa\leq 0$. We choose $r<d(o,p)/4$.
For $u,v\in B_r(p)$ we have by the euclidean Law of Sines
$\sin(\angle_o(u,v))\leq\frac{d(u,v)}{d(o,v)}<\frac{d(u,v)}{d(o,p)-r}$
and $\angle_o(u,v)<\pi/2$.
For $\kappa>0$ we choose $r<D_\kappa/4$
and we obtain similarly that
$\sin(\angle_o(u,v))\leq\frac{\sin(\sqrt\kappa d(u,v))}{\sin(\sqrt\kappa d(o,v))}
<\frac{\sin(\sqrt\kappa d(u,v))}{\sin(\sqrt\kappa(d(o,p)-r))}$
and $\angle_o(u,v)<\pi/2$.
\qed
\end{Lem}
\begin{Num}\textbf{The blow-up at $o$.}\label{BlowUp}\\
Let $(\bar X,d)$ denote the metric completion of $(Y,d)$ (equivalently,
of $(X,d)$). We denote the completion of $(Y,d_o)$ by $(\hat X,d_o)$ and
call it the \emph{blow-up}
\footnote{This is not the blow-up in the sense of
Lytchak \cite{Lyt}.} of $X$ at $o$. From the
the $1$-Lipschitz maps $(Y,d_o)\rTo(Y,d)$ and 
$(Y,d_o)\rTo(Y,\angle_o)$ we obtain $1$-Lipschitz maps
\[
q:(\hat X,d_o)\rTo(\bar X,d)\text{ and }\hat\rho:(\hat X,d_o)\rTo(\Sigma_oX,\angle_o).
\]
We note also that $d$ and $\angle_o$ extend to continuous pseudometrics on
$\hat X$, which we denote by the same symbols.
\end{Num}
\begin{Prop}
The preimage $q^{-1}(o)$ is the space of directions
$\Sigma_oX$. The restriction $q:\hat X\setminus q^{-1}(o)\rTo\bar X\setminus o$
is a homeomorphism. Thus we have
\[
\hat X=(\bar X\setminus o)\,\dot\cup\, \Sigma_oX
\]
and $\Sigma_oX\rInto\hat X$ is a closed isometric embedding.

\proof
Let $(x_n)_{n\in\mathbb N}$ and $(y_n)_{n\in\mathbb N}$
be Cauchy sequences in $(Y,d_o)$
with $\lim_nd(x_n,o)=\lim_nd(y_n,o)=0$.
Then $\lim d_o(x_n,y_n)=\angle_o(x_n,y_n)$. Therefore their
limits $x$ and $y$ in $\hat X$ have distance $d_o(x,y)=\angle_o(x,y)$.
Both sequences are also Cauchy with respect to $\angle_o$. So
$q^{-1}(o)$ is a complete subset of $\Sigma_oX$.
If $y\in Y$ represents a point in $\Sigma_oX$, then
the $d_o$-Cauchy sequence $(\frac1ny)_{n\in\mathbb N}$ represents
the same point. Thus $q^{-1}(o)$ is dense in $\Sigma_oX$.
\qed
\end{Prop}
\begin{Lem}
For $s\in[0,1]$ and $x\in\hat X$ we put 
\[
\hat h_s(x)=\begin{cases}
        sx& \text{for }x\in \bar X\setminus o\text{ and }s>0\\ 
        \hat\rho(x)& \text{else.}
       \end{cases}
\]
Then $(x,s)\mapstoo \hat h_s(x)$ is  a $1$-Lipschitz 
strong deformation retraction of $\hat X$ onto $\Sigma_oX$.

\proof
We have $\angle_o(u,v)=\angle_o(\hat h_s(u),\hat h_{s'}(v))$ for all
$u,v\in\hat X$ and $s,s'\in[0,1]$.
From the comparison triangle in $M_\kappa$ we see that
$d(su,s'v)^2\leq d(u,v)^2+(s-s')^2$, whence
$d_o(su,s'v)^2\leq d_o(u,v)^2+(s-s')^2$.
Thus $(x,s)\mapstoo \hat h_s(x)$ is $1$-Lipschitz.
\qed
\end{Lem}
We remark that this homotopy can be restricted to the subset
$Y\cup\Sigma_oX\subseteq\hat X$.

\section{The homotopy type of the space of directions}

In this section we prove that the space of directions at a point $o$
in a CAT$(\kappa)$ space $X$
has the homotopy type of a small punctured neighborhood of $o$.
We first collect some facts from homotopy theory.
Unless specified otherwise, all simplicial complexes are endowed with
the weak topology (the CW topology). This allows us to define
continuous maps simplex-wise. We start with a fact from
obstruction theory.
\begin{Lem}[Acyclic carriers]
\label{Acyclic}
Let $K$ be a simplicial complex and $Y$ a topological space.
Let $f:K^{(0)}\rTo Y$ be a map from the $0$-skeleton (the
vertex set) of $K$ to $Y$. Suppose that there is a map $C$ that
assigns to every simplex $A\subseteq K$ a subset $C_A\subseteq Y$
with the following properties.
\begin{enumerate}
\item[(1)] If $A\subseteq K$ is a simplex, then $C_A$ contains $f(A\cap K^{(0)})$.
\item[(2)] If $A\subseteq B$ are simplices, then $C_A\subseteq C_B$.
\item[(3)] $\pi_*(C_A,p)=0$ for all $p\in C_A$.
\end{enumerate}
Then $f$ has a continuous extension $f:K\rTo Y$ such that
$f(A)\subseteq C_A$ holds for all simplices $A\subseteq K$.

\proof
The map $f$ is defined inductively on the $m$-skeleton.
Suppose that $f:K^{(m-1)}\rTo Y$ is already defined.
Let $A\subseteq K$ be an $m$-simplex and put 
$\partial A=K^{(m-1)}\cap A$. Choose $p\in f(\partial A)\subseteq C_A$.
Since $\pi_{m-1}(C_A,p)=0$, there exist an extension of the map 
$f|_{\partial A}:\partial A\rTo C_A\subseteq Y$
over $A$; see \cite[VI.6.5]{HuHom} or \cite[V.5.14]{Whi}.
These maps fit together to a continuous map $f:K^{(m)}\rTo Y$.
\qed
\end{Lem}
We call a function $C$ satisfying (1), (2), (3) in the previous lemma
a \emph{topological acyclic carrier} for $f$. Our first application
is as follows~%
\footnote{After the present paper was completed, Urs Lang informed
me that \cite[Thm.~7]{Girolo} contains this result for complete
CAT$(\kappa)$ spaces. See also the remarks in \cite[4.2]{Fi}
and \cite[I.7A.15,II.5.13]{BH}.
The result is also proved in \cite[Lem~1.1]{Ontaneda} and in
\cite[Sec.~6]{LyS}---essentially by the same characterization of ANRs.}.
\begin{Thm}
\label{CATImpliesANR}
Let $Z$ be a metric space with curvature bounded above. Then every
open subset of $Z$ is an ANR (an absolute neighborhood retract).

\proof
It suffices to show that $Z$ is locally an AR (an absolute retract); see \cite[III.8.1]{Hu}.
Let $o\in Z$ and choose $R>0$ in such a way that $\bar B_R(o)$
is a CAT$(\kappa)$ space of radius less than $D_\kappa/2$.
We show that $X$ is an AR (absolute retract).
Let $K$ denote the simplicial complex on all finite subsets of
$X$. Suppose that $A\subseteq K$ is a simplex spanned by the points
$a_0,\ldots,a_m\in X$.
Let $C_A$ denote the convex hull of  $a_0,\ldots,a_m$ in $X$. Then
$C_A$ is contractible. In this way we obtain a topological acyclic carrier
for the identity map $X=K^{(0)}\rTo X$.
By Lemma~\ref{Acyclic}, there exists a continuous map $q:K\rTo Z$ such that
the image of every simplex $A\subseteq K$ is contained in $C_A$.
Since the radius of $C_A$ is bounded by the diameter of the set
$\{a_0,\ldots,a_m\}$, the map $q$ has the following two properties:
\begin{enumerate}
\item[(1)] $q$ is surjective.
\item[(2)] If $(A_j)_{j\in\mathbb N}$ is a sequence of simplices in $K$ such that
  $(q(A_j\cap K^{(0)}))_{j\in\mathbb N}$ converges to $x\in X$, then 
  $(q(A_j))_{j\in\mathbb N}$ also converges to $x$.
\end{enumerate}
By  Wojdys\l awski's Theorem~\ref{Woj}, $X$ is an AR.
\qed
\end{Thm}
Next, we recall the following version of the Whitehead Theorem.
\begin{Thm}[Whitehead]\label{Whitehead}
Let $f:X\rTo Y$ be a continuous map between topological spaces.
Then $f$ is a weak equivalence if and only if for every
finite simplicial complex $K$ the induced map
\[
f_*:[K,X]{\rTo}[K,Y]
\]
between free homotopy sets is bijective.

\proof
The bijectivity is sufficient by \cite[IV.7.16]{Whi}.
There, the result is stated under
the stronger assumption that bijectivity
holds for every CW complex $K$. However, the only spaces $K$
used in the proof in \textit{loc.cit.}~are spheres.
The bijectivity is a necessary condition by \cite[IV.7.17]{Whi}.
\qed
\end{Thm}
Armed with \ref{CATImpliesANR} and Whitehead's Theorem, we now proceed to
the proof of Theorem~A. The following Lifting Lemma is the crucial step.
Related results can be found in \cite{Kl}. However, our topological
approach avoids
all completeness assumptions (for example, we do not use barycentric
simplices).
\begin{Lem}[The Lifting Lemma]
\label{LiftingLemma}
Let $X$ be a CAT$(\kappa)$ space and suppose that $o\in X$ is
a point with $X\subseteq\bar B_{D_\kappa/2}(o)$.
Put $Y=X\setminus o$ and let $\rho:Y\rTo\Sigma_oX$ be as in \ref{SigmaX}.
Assume that $K$ is a finite simplicial complex and that $g:K\rTo\Sigma_oX$
is a continuous map. Then there exists for every $\eps>0$ a 
continuous map $f:K\rTo Y$ such that $\rho\circ f$ is
$\eps$-near to $g$.

\proof
We may assume that $\kappa>0$, since we are interested in a local question.
We may also assume that $\eps<\pi/2$. We choose $\alpha>0$ such that
$\alpha <\max\{1/2,\sin(\eps/4)\}$. We subdivide
$K$ in such a way that the $g$-image of every simplex $A\subseteq K$
has diameter less than $\alpha/3$. Next, we choose for every vertex
$a$ of $K$ a point $y_a\in Y$ such that $\angle_o(\rho(y_a),g(a))<\alpha/3$
(this is possible because $\rho(Y)$ is dense in $\Sigma_oX$).
The vertex set $K^{(0)}$ is finite, so moving points along geodesic
segments towards $o$,
we can arrange that all $y_a$ have the same positive distance $r$ from $o$.
If the vertices $a,b$ belong to the same simplex $A$ in $K$, then 
\[
\angle_o(y_a,y_b)<\angle_0(y_a,g(a))+\angle_o(g(a),g(b))+\angle_o(g(b),y_b)<\alpha.
\] 
Since
\[
\sin(\angle_o(y_a,y_b)/2)=\lim_{s\to0}{d(sy_a,sy_b)}/{(2rs)}<\sin(\alpha/2)<\alpha/2,
\]
we have $d(sy_a,sy_b)<sr\alpha$ for all sufficiently small $s>0$.
Since $K^{(0)}$ is finite, there exists an $s_0>0$ such that
the following holds for all $s\in(0,s_0)$:
\begin{enumerate}
 \item[(*)] If $\{a_0,\ldots,a_m\}$ span a simplex $A$ in $K$, then the
convex hull $C_A$ of $\{sy_{a_0},\cdots,sy_{a_m}\}$ has radius less than
$sr\alpha$. 
\end{enumerate}
Since $\alpha<1/2$ and $d(sy_{a_j},o)=sr$, this implies in particular that $C_A$ does not contain $o$, 
i.e. $C_A\subseteq Y$. We put $f(a)=sy_a$. In this way we have constructed a topological
acyclic carrier and we may apply Lemma \ref{Acyclic}. We obtain a continuous extension
$f:K\rTo Y$ , such that $f(A)\subseteq C_A$ for every simplex $A\subseteq K$.
Such an extension exists for every $s<s_0$. In particular, we may assume that
$\sqrt\kappa sr\alpha<\pi/2$.
If $z\in Y$ is any point with $d(z,sy_a)< sr\alpha$, then we have by
the Law of Sines \ref{LawOfSines}
\[
\sin(\angle_o(z,sy_a))<\frac{\sin(\sqrt\kappa sr\alpha)}{\sin(\sqrt\kappa sr)}.
\]
Since $\alpha<\sin(\eps/4)$, we see (using l'H\^opital's rule) that for all
sufficiently small $s$, we have 
$\frac{\sin(\sqrt\kappa sr\alpha)}{\sin(\sqrt\kappa sr)}<\sin(\eps/4)$.
Thus we can choose $s$ in such a way that $\angle_o(z,sy_a)<\eps/4$
holds for all vertices $a$ and all $z\in B_{sr\alpha}(sy_a)$. 
For such an $s$ we consider the map $f:K\rTo Y$ constructed before.
If $x$ is in the simplex $A\subseteq K$ and $a\in A$ is a vertex, then 
\begin{align*}
\angle_o(\rho(f(x)),g(x))&\leq
\angle_o(\rho(f(x)),\rho(f(a)))+
\angle_o(\rho(f(a)),g(a))+
\angle_o(g(a),g(x))\\
&<\eps/4+\alpha/3+\alpha/3+\alpha/3\\
&<\eps.
\end{align*}
\qed
\end{Lem}
\begin{Thm}
\label{TheoremA}
Let $X$ be a CAT$(\kappa)$ space and suppose that $o\in X$ is
a point with $X\subseteq\bar B_{D_\kappa/2}(o)$.
Put $Y=X\setminus o$ and let $\rho:Y\rTo\Sigma_oX$ be as in \ref{SigmaX}.
Then $\rho$ is a homotopy equivalence.

\proof
Let $K$ be a finite simplicial complex. We show that 
$\rho_*:[K,Y]\cong[K,\Sigma_oX]$ is a bijection. Suppose that $g:K\rTo\Sigma_oX$ is
continuous. By the Lifting Lemma
\ref{LiftingLemma} there is a continuous map $f:K\rTo Y$
such that $\rho\circ f$ is $\pi/2$-near to $g$. Therefore
$\rho\circ f$ is homotopic to $g$ by the geodesic homotopy in $\Sigma_oX$.
This shows that $\rho_*$ is surjective. 

Now suppose that $h_0,h_1:K\pile{\rTo\\ \rTo}Y$
are continuous and that $\rho\circ h_0$ and $\rho\circ h_1$ are
homotopic in $\Sigma_oX$. Thus there is a homotopy $g:K\times[0,1]\rTo\Sigma_oX$ with
$\rho\circ h_i=g_i$, for $i=0,1$. 
By the Lifting Lemma
\ref{LiftingLemma}, there exists a continuous map $f:K\times[0,1]\rTo Y$ such that
$\rho\circ f$ is $\pi/2$-near to $g$. For $x\in K$ and $i=0,1$ 
we have 
\begin{align*}
\angle_o(h_i(x),f_i(x))&=\angle_o(\rho(h_i(x)),\rho(f_i(x)))\\
&=\angle_o(g_i(x),\rho(f_i(x)))\\
&\leq\pi/2.
\end{align*}
In particular, $o$ is not in the geodesic segment
$[h_i(x),f_i(x)]\subseteq X$, since otherwise we would have
$\angle_o(h_i(x),f_i(x))=\pi$. Thus the geodesic homotopy in $X$ between
$h_i$ and $f_i$ takes its values in $Y$. Therefore
$h_0\simeq f_0\simeq f_1\simeq h_1$ in $Y$. This shows that $\rho_*$
is injective. By Whitehead's Theorem~\ref{Whitehead}, $\rho$ is a weak homotopy equivalence.
Since $Y$ and $\Sigma_oX$ are ANRs by \ref{CATImpliesANR}, both spaces
have the homotopy types of CW complexes; see \ref{Mard}. A weak homotopy equivalence
between such spaces is a homotopy equivalence by Whitehead's (other) Theorem; 
see \cite[V.3.5]{Whi}.
\qed
\end{Thm}
This finishes the proof of Theorem~A in the introduction
\footnote{Bruce Kleiner informed me that he had obtained the following result
some time ago (unpublished):
for any open set $W\subseteq \Sigma_oX$ and $\eps>0$ small enough, the map 
$(B_\eps(o)\setminus o)\cap\rho^{-1}(W)\rTo W$
is a homotopy equivalence. In fact, our proof above could easily be modified
in order to obtain this stronger result.}.

\section{The local homology sheaf of a space}

Let $X$ be a topological space. Recall that a \emph{presheaf} on $X$ is a
cofunctor
\footnote{With values in a some given category.}
on the category of open sets of $X$. An important example is
the map which assigns to
an open set $U\subseteq X$ the relative (singular)
homology group $H_*(X,X\setminus U)$.
The sheaf generated by this presheaf is the
\textit{local homology sheaf} $\mathcalligra H\,_*(X)$ of $X$. Its stalk at $p\in X$ is the
local singular homology
\[
\mathcalligra H\,_*(X)_p=\lim_{U\ni p}H_*(X,X\setminus U)=H_*(X,X\setminus p);
\]
see \cite[p.~7]{Bre}.
Every relative $k$-cycle $\sigma\in H_k(X,X\setminus U)$
induces a section $p\mapsto \sigma_p\in \mathcalligra H_k\,(X)_p$ over $U$ via the restriction map
$H_k(X,X\setminus U)\rTo H_k(X,X\setminus p)$. The \textit{support}
of $\sigma$ is 
\[
\supp(\sigma)=\{p\in U\mid\sigma_p\neq 0\}.
\]
This set is closed in $U$ (by general facts about sheaves, or more directly
because singular homology satisfies the axiom of compact supports).
We note also that $\supp(\sigma)$ is contained in the image of
any relative singular cycle representing $\sigma$.

On $X$ we have also the set-valued presheaf
\[
Closed:U\mapstoo\{A\cap U\mid A\subseteq X\text{ is closed}\}.
\]
The sheaf $\mathcalligra{Closed\,}\ $ generated by $Closed$ has as its stalk 
$\mathcalligra{Closed\,}_p$
the set of germs of closed sets near $p$.
This stalk $\mathcalligra{Closed\,}_p$ is a poset and a distributive lattice
\footnote{It is a 'pointless topological space'.};
see~\cite[I\,\S6.9,I\,\S6.10]{Bki} or \cite{Izumi}.
It is clear that $\supp$ is a natural transformation of
presheaves 
\[
\supp:H_*(X,X\setminus U)\rTo Closed(U).
\]
and sheaves, in particular we have an induced map on the stalks
\[
\supp:\mathcalligra H\,_*(X)_p
\rTo \mathcalligra{Closed\,}_p\,.
\]

\section{Cycles in spherical buildings}

We start with some elementary observations about simplicial complexes.
(We recall our convention that we endow simplicial complexes with the
weak topology.)
Let $K$ be an $m$-dimensional simplicial complex and let
$\sigma$ be an $m$-cycle in the homology $H_m(K)$. 
We are intersted in the support of $\sigma$ over $K$.
\begin{Lem}\label{SupportIsPure}
The support $\supp(\sigma)\subseteq K$ is a pure
\footnote{A complex is pure if all maximal simplices have the same dimension.}
$m$-dimensional subcomplex.

\proof
We represent $\sigma$ as a finite sum $\beta_1C_1+\cdots\beta_kC_k$, where the $C_j$ are
distinct $m$-simplices in $K$ and the $\beta_j$ are nonzero coefficients. For each
interior point $p\in C_j$, the cycle $\sigma$ maps in 
$H_m(C_j,C_j\setminus p)\cong H_m(K,K\setminus p)\cong\ZZ$
to $\beta_j$ times a generator (by excision). Thus $C_1\cup\cdots\cup C_k\subseteq \supp(\sigma)$.
On the other hand, the support set cannot be bigger than the union of the simplices
representing $\sigma$.
\qed
\end{Lem}
The \emph{join} of two simplicial complexes $K$, $L$ is denoted $K*L$.
The complex $K*\SS^0$ is the same as the unreduced suspension of $K$.
\begin{Lem}
Under the Mayer-Vietoris isomorphism $H_{m+1}(K*\SS^0)\rTo H_m(K)$,
the preimage of $\sigma\in H_m(K)$ has $\supp(\sigma)*\SS^0$
as its support set.

\proof
Let $\SS^0=\{u,v\}$.
As in the previous proof we represent $\sigma$ as a finite sum
$b=\beta_1C_1+\cdots\beta_kC_k$, where the $C_j$ are
distinct $m$-simplices and the $\beta_j$ are nonzero coefficients.
Consider the two $m+1$-chains $u*b=\beta_1(u*C_1)+\cdots+\beta_k(u*C_k)$
and $b*v=\beta_1(C_1*v)+\cdots+\beta_k(C_k*v)$. Then
$\partial(u*b)=-b$ and $\partial(b*v)=b$, so
$u*b+b*v$ is a cycle in $K*S^0$, with support
$\supp(\sigma)*\SS^0$. We claim that it maps to $b$.
Indeed, the map $H_{m+1}(K*\SS^0)\rTo H_m(K)$ is the composite
\[
H_{m+1}(K*\SS^0)\rTo^\cong H_{m+1}(K*\SS^0,K*u)\lTo^\cong H_{m+1}(K*v,K)\rTo^\cong H_m(K);
\]
see~\cite[\S15]{ES}.
Tracing the cycle $u*b+b*v$ along these arrows, we end up at $b$.
\qed
\end{Lem}
\begin{Cor}\label{SupportInJoin}
The preimage of $\sigma\in H_m(K)$ under the isomorphism $H_{m+n+1}(K*\SS^n)\rTo H_m(K)$
has as its support $\supp(\sigma)*\SS^n$.

\proof
This is an iterated application of the previous lemma,
since $\SS^0*\SS^k\cong\SS^{k+1}$.
\qed
\end{Cor}
\begin{Num}\textbf{Spherical buildings}\label{SphericalBuildings}\\
A \emph{spherical building} $\Delta$ is a simplicial complex together with
a collection $Apt(\Delta)$ of subcomplexes, called \emph{apartments}, which are
isomorphic to a fixed spherical Coxeter complex $\Sigma$ (a certain triangulated
sphere). The apartments have to satisfy the following compatibility conditions.
\begin{enumerate}
\item[(B1)]
For any two simplices $C,D\subseteq \Delta$, there is an apartment $A\subseteq\Delta$
containing $C,D$.
\item[(B2)]
If $A,A'\subseteq \Delta$ are apartments containing the simplices $C,D$, then there is
a combinatorial (type preserving) isomorphism $A\rTo A'$ fixing $C$ and $D$.
\end{enumerate}
The maximal simplices are also called \emph{chambers}.
If every non-maximal simplex is contained in at least three apartments, then
the building is called \emph{thick}. Buildings in the sense of the present
definition are sometimes called \emph{weak buildings}.
We refer to \cite{TitsLNM} and to
\cite{AB,TitsComo,WeSph} for details and to \cite{CL} for a nice characterization
from a metric viewpoint.
\end{Num}
Let now $\Delta$ be a (possibly weak)
spherical building of dimension $m$ (viewed as a simplicial complex
and endowed with the CW topology). By the Solomon-Tits Theorem,
$\Delta$ has the homotopy type of a sum of $m$-spheres. This can be
made more precise as follows. Let $C_0$ be a chamber (an $m$-simplex) and let
$M$ denote the collection of all chambers opposite $C_0$. Each chamber
$C\in M$ determines together with $C_0$ a unique apartment $A\simeq\SS^m$.
The complex $N=\Delta-\bigcup\{\mathrm{int}(C)\mid C\in M\}$ consisting
of all simplices which are not in $M$ is contractible
and thus 
\[
\Delta\simeq\Delta/N\cong\bigvee_M\SS^m.
\]
See \cite[4.12]{AB} for details.
We now inspect the cycles in the top-dimensional singular homology of $\Delta$.
For each apartment
$A\subseteq\Delta$ containing $C_0$ we have the cycle $[A]\in H_m(\Delta)$ 
given by the fundamental class
\footnote{We fix an orientation of the simplex $C_0$ and choose the orientation of
the $m$-sphere $A$ accordingly.} of $A$. These cycles $[A]$ form a basis for the
free $\ZZ$-module $H_m(\Delta)$ and the support of
$[A]$ is $A$, since the fundamental class of $A$ consists of all simplices in $A$.
\begin{Lem}
Let $\Delta$ be a thick spherical $m$-dimensional building and
let $n\geq -1$. Put $\SS^{-1}=\emptyset$.
Let $\mathcal S$ denote the lattice consisting of all finite intersections
and unions
of support sets of cycles in $H_{m+n+1}(\Delta*\SS^n)$.
Then $\mathcal S$ has a unique minimal element, $\emptyset*\SS^n$.
If $a$ is a simplex in $\Delta$, then $a*\SS^n$ is in $\mathcal S$.

\proof
By \ref{SupportInJoin}, every support set of a cycle in 
$H_{m+n+1}(\Delta*\SS^n)$ is of the form
$\supp(\sigma)*\SS^n$, for some cycle $\sigma\in H_m(\Delta)$.
Therefore the first claim follows. 
Since $\Delta$ is thick, every simplex $a\in\Delta$ can be
written as the intersection of finitely many apartments in $\Delta$.
This shows the second claim.
\qed
\end{Lem}
\begin{Num}\label{RecoverBuilding}
We can recover the combinatorial structure of $\Delta*\SS^n$
from $\mathcal S$ as follows. 
Let us call an element of a poset lattice \emph{indecomposable} if it is not
a union of finitely many strictly smaller elements.
All elements of $\mathcal S$ are of the form $K*\SS^n$, where $K$ is a finite subcomplex of
$\Delta$. We have seen above that every simplex $a$ of $\Delta$ occurs
in this way in $\mathcal S$ as $a*\SS^n$. Obviously,
these sets $a*\SS^n$ are precisely the indecomposable elements of $\mathcal S$.
The combinatorial structure of the underlying simplicial complex of a thick
spherical building determines the buildings uniquely: there is precisely
one set of apartments which turns $\Delta$ into a spherical building; 
see \cite[3.15]{TitsLNM}.
\end{Num}
\begin{Num}\textbf{Thick reductions}\label{ThickReduction}\\
The reason why we consider joins of spherical buildings and spheres is as
follows.
If $\Delta$ is a weak spherical building, then there exists a thick
spherical building $\Delta_0$ 
\footnote{Strictly speaking, we have to allow here that $\Delta_0=\emptyset=\SS^{-1}$.}
and a number $n\geq -1$ such that
a simplicial subdivision of $\Delta$ is isomorphic to $\Delta_0*\SS^n$;
see \cite[1.3]{Cap}, \cite[3.8. or 7.1]{CL}, \cite[3.7]{KL}, \cite{Scha}.
\end{Num}
\begin{Num}
So far, we have considered spherical buildings as simplicial complexes
endowed with the weak (simplicial) topology. A spherical building carries
also another topology coming from its natural CAT$(1)$ metric, 
see~\cite[II.10A.4]{BH}. The identity is a continuous  map from the
weak topology to the metric topology. If $\Delta$ is infinite, then this
is not a homeomorphism. However, it is always a homotopy equivalence.
For a proof, see
Section I.7 in \cite{BH}, in particular \cite[Exc.~I.7A.11(2)]{BH},
combined with \cite[Thm.~1]{Dow} or \cite[5.4.6]{LW}. The metric topology yields
the same cycles and support sets in $\Delta$ as the weak topology.
That the support sets remain unchanged can be seen as follows:
the metric topology and the CW topology coincide in the interior of
every chamber.
\end{Num}

\section{Topological rigidity of euclidean buildings}

We recall briefly the definition of a euclidean building. We refer to
\cite{TitsComo} and to \cite{KL,KrWe,Par} for details and further
results.
\begin{Num}\textbf{Euclidean buildings}\label{EBAxioms}\\
Let $W$ be a spherical Coxeter group and $W\RR^n$ the corresponding
affine Weyl group (we extend $W$ by the full translation group $(\RR^n,+)$).
From the reflection hyperplanes of $W$ we
obtain a decomposition of $\RR^n$ into \emph{walls}, \emph{half spaces},
\emph{Weyl chambers} (a Weyl chamber
is a fundamental domain for $W$---these are Tits' \emph{chambres vectorielles}) and 
\emph{Weyl simplices} (Tits' \emph{facettes vectorielles}).

Let $X$ be a metric space. A \emph{chart} is an isometric
embedding $\phi:\mathbb E^n\rTo X$, and its image is called an \emph{affine apartment}.
We call two charts $\phi,\psi$
\emph{$W$-compatible} if $Y=\phi^{-1}(\psi(\mathbb E^n))$ is convex 
(in the Euclidean sense) and
if there is an element $w\in W\RR^n$ such that
$\psi\circ w|_Y=\phi|_Y$ (this condition is void if $Y=\emptyset$). 
We call a metric space $X$ together with a collection $\mathcal A$ of
charts a \emph{Euclidean building} if it has the following
properties.
\begin{enumerate}
\item[(A1)] For all $\phi\in\mathcal A$ and $w\in W\RR^n$, the composition
$\phi\circ w$ is in $\mathcal A$.
\item[(A2)] The charts are $W$-compatible.
\item[(A3)] Any two points $x,y\in X$ are contained in some affine apartment.
\end{enumerate}
The charts allow us to map Weyl chambers, walls and half spaces into $X$;
their images are also called Weyl chambers, walls and half spaces. The first
three axioms guarantee that these notions are coordinate independent.
\begin{enumerate}
\item[(A4)] If $C,D\subseteq X$ are Weyl chambers, then there is
an affine apartment $A$ such that the intersections
$A\cap C$ and $A\cap D$ contain Weyl chambers.
\item[(A5')] For every apartment $A\subseteq X$ and every $p\in A$ there is
a $1$-Lipschitz retraction $h:X\rTo A$ with $h^{-1}(p)=\{p\}$.
\end{enumerate}
Condition (A5') may be replaced by the following condition:
\begin{enumerate}
\item[(A5)] If $A_1,A_2,A_3$ are affine apartments which intersect pairwise
in half spaces, then $A_1\cap A_2\cap A_3\neq\emptyset$.
\end{enumerate}
See \cite{Par} for a thorough discussion of different sets of axioms.
\end{Num}
We assume now that $X$ is a euclidean building with $n$-dimensional apartments.
\begin{Lem}\label{SplitInjection}
Let $A\subseteq X$ be an apartment and
let $p\in A$. Then $H_n(A,A\setminus p)\rTo H_n(X,X\setminus p)$ is a split injection.

\proof
This is clear from Axiom (A5').
\qed
\end{Lem}
\begin{Cor}
\label{RestrictionCor}
Let $r>0$ and let $\sigma$ be a generator of $H_n(A,A\setminus B_r(p))\cong\ZZ$.
The support of the image of $\sigma$ in $H_n(X,X\setminus B_r(p))$ is
$A\cap B_r(p)$.

\proof
We represent $\sigma$ as the fundamental class of the 
closed $n$-disk $\bar B_r(p)\cap A$ relative to its boundary, an $n-1$-sphere.
This class generates
$H_n(A,A\setminus B_r(p))\cong\ZZ$.
We conclude that $\supp(\sigma)\subseteq A\cap B_r(p)$. For every $o\in A\cap B_r(p)$
we have a commutative diagram
\begin{diagram}[size=2em]
H_n(A,A\setminus B_r(p)) &\rTo& H_n(X,X\setminus B_r(p))\\
\dTo^\cong &&\dTo\\
H_n(A,A\setminus o) &\rTo_{\text{\small injective}}& H_n(X,X\setminus o).
\end{diagram}
Thus $o$ is in the support.
\qed
\end{Cor}
\begin{Num}
Let $o\in X$. We recall the following facts (essentially, all due to Tits); 
see \cite[2.10]{Par}
\cite[4.2.2 and 4.4.]{KL} or \cite[1.14 and p.~20]{Par} for proofs.
\begin{enumerate}
 \item[(a)] A euclidean building is a (not necessarily complete) CAT$(0)$ space.
 \item[(b)] The space of direction $\Sigma_oX$ at any point $o\in X$
  is the CAT$(1)$ realization of an $n-1$-dimensional
 (possibly weak) spherical building $\Delta=\Delta_o$ of type $(W_0,I)$.
 \item[(c)] Let $\mathcal W_o$ denote the collection of all $o$-based Weyl
simplices in $X$. The image of $\mathcal W_o$ in $\mathcalligra{Closed\,}_p$
is (as a poset) precisely the simplicial complex $\Delta_o$.
\end{enumerate}
\end{Num}
The next result is essentially \cite[6.2.3]{KL}.
\begin{Prop}
\label{SupportGerm}
Let $\sigma\in H_n(X,X\setminus B_r(o))$. Then there exists $\eps\in (0,r)$
such that the following holds. There are $o$-based Weyl chambers
$C_1,\ldots,C_k$ such that the support of $\sigma$ in $B_\eps(o)$
is $\supp(\sigma)=B_\eps(o)\cap(C_1\cup\cdots\cup C_k)$.

\proof
Since $X$ is contractible, we have by \ref{TheoremA}
\[
H_n(X,X\setminus B_r(o))\rTo
H_n(X,X\setminus o)\rTo_\cong^\partial H_{n-1}(X\setminus o)\rTo_\cong^{\rho_*} H_{n-1}(\Sigma_oX).
\]
Let $K\subseteq \Sigma_oX$ denote the support of the image of $\sigma$ in $\Sigma_oX$.
By \ref{SupportIsPure} and (b) above, $K$ is a finite pure subcomplex of the
spherical building $\Delta_o$.
We fix a finite set $\mathcal K_o\subseteq\mathcal W_o$ of $o$-based
Weyl simplices representing the simplices in $K$. Since $K$ is finite,
we can find by (c) a number $s\in (0,r)$ such that for all
simplices $A',B'\subseteq K$ with corresponding Weyl simplices
$A,B\in\mathcal K_o$, we have 
\[
A'\subseteq B'\text{ if and only if }
\bar B_s(o)\cap A\subseteq \bar B_s(o)\cap B.
\]
Thus $K_s=\{x\in \bigcup\mathcal K_o\mid d(x,o)=s\}$ is a simplicial
complex which is isomorphic to $K$ under the canonical map
$\rho:X\setminus o\rTo\Sigma_oX$. The set
$K_{\leq s}=\{x\in \bigcup\mathcal K_o\mid d(x,o)\leq s\}$ is a finite
$n$-dimensional simplicial complex isomorphic to the cone over $K_s$.

We now write the image of $\sigma$ in $H_{n-1}(\Sigma_oX)$
as a linear combination of distinct chambers $b=\beta_1C_1'+\cdots+\beta_kC_k'$
with nonzero coefficients $\beta_j$. Let $C_j\in\mathcal K_o$ denote the
Weyl chamber representing $C_j'$, put $\tilde C_j=C_j\cap\bar B_s(o)$
and consider the relative cycle
$\beta_1\tilde C_1+\cdots+\beta_k\tilde C_k$ in $H_m(K_{\leq s},K_s)$.
Let $\tau$ denote its image in $H_m(X,X\setminus B_s(o))$.
If $p$ is an interior point in the $m$-simplex
$\tilde C_j$, then $\tau$ restricts by excision and \ref{RestrictionCor}
to $\beta_j$ times a generator
of $H_m(\tilde C_j,\tilde C_j\setminus p)$. 
Thus $\supp(\tau)=K_{\leq s}\cap B_s(o)$. Moreover,
$\tau$ has (by construction) the same image in $H_{n-1}(\Sigma_oX)$ as
$\sigma$. Therefore there exists
$\eps\in(0,s)$ such that $\tau$ and $\sigma$ have the same support set
in $B_\eps(o)$ (because the complements of supports are open).
\qed
\end{Prop}
\begin{Cor}
\label{RecoverFromGerms}
Consider the image of $\supp:\mathcalligra H_*\,(X)_o
\rTo \mathcalligra{Closed\,}_o$. Let $\mathcal S_o\subseteq \mathcalligra{Closed\,}_o$
denote the lattice consisting of all finite intersections and unions of germs of support sets.
The indecomposable elements in $\mathcal S_o$ form the simplices of the
thick part of the reduction of
$\Delta_o$, see \ref{RecoverBuilding}. The germ of a support set comes from the fundamental
class of an apartment $A$ containing $o$ if and only if it can be written as
a union of indecomposables representing an apartment in the thick reduction of $\Delta_o$.
\end{Cor}
In particular, we obtain the following result.
\begin{Thm}
Let $f:X\rTo X'$ be a homeomorphism between two euclidean buildings
with $n$ and $n'$-dimensional apartments, respectively. Then $n=n'$. If 
$A\subseteq X$ is an apartment and $o\in A$, then there exists
a small neighborhood $U\subseteq A$ of $o$ so that $f(U)$ is contained in
an apartment $A'$ of $X'$.

\proof
By \ref{SplitInjection}, the number $n$ is characterized by the
fact that $H_j(X,X\setminus o)=0$ for $j\neq n$ (we can of course also
use the equality $n=\dim(X)=\dim(X')=n'$ established in \ref{DimensionThm}).
Thus $n=n'$.
From \ref{RecoverFromGerms} we see that $\Delta_o$ is a sub-poset of
$\mathcalligra{Closed\,}_o$ that can be read off from
$\mathcalligra{H}\,_n(X)_o$. In particular, we can read off
which cycles in $\mathcalligra{H}\,_n(X)_o$ arise from
fundamental classes of apartments. The claim follows from \ref{SupportGerm}.
\qed
\end{Thm}
An immediate corollary is the following result due to Kleiner and Leeb \cite[6.4.2]{KL}.
\begin{Cor}[Kleiner-Leeb]
\label{TopologicalRigidity}
Let $X,X'$ be euclidean buildings and suppose that
$f:X\rTo X'$ is a homeomorphism. Let $A\subseteq X$ be an
apartment. If $X'$ is complete, then $f(A)$ is an apartment in the maximal
atlas $\mathcal A_{\max}$ of $X'$.

\proof
For each $o\in A$ there is a small neighborhood $U\subseteq A$
of $o$ that is contained in some apartment $A'\subseteq X'$. Thus
$f(A)$ is a complete simply connected metric space that
is locally isometric to euclidean space $\mathbb E^n$.
In particular, $f(A)$ is a flat complete simply connected
Riemannian manifold. Such a manifold is isometric to
euclidean space; see \cite[V.4.1]{Kob}. It follows that
$f(A)$ is an apartment in $\mathcal A_{\max}$; see
\cite[4.6]{KL} or \cite[2.25]{Par}.
\qed
\end{Cor}
The following example shows that the completeness assumption in \ref{TopologicalRigidity}
is, apparently, not superfluous. We also remark that there exist higher
dimensional non-complete Bruhat-Tits buildings;
some algebraic conditions are discussed in \cite[7.5]{BruTi}.
\begin{Num}\textbf{Example}\label{TreeExample}
Consider the $\RR$-trees $T_1=(-1,1)\times\RR$ and $T_2=\RR\times\RR$, both endowed with the
metric 
\[
d((x,y),(x',y'))=\begin{cases}
                 |y-y'|&\text{if }x=x'\\
                 |y|+|x-x'|+|y'|&\text{if }x\neq x'.
                 \end{cases}
\]
Both trees are $1$-dimensional euclidean buildings and $T_2$ is complete
(but $T_1$ is not complete).
Stretching the $x$-axis, we obtain
a homeomorphism $T_2\rTo T_1$. In $T_2$, the $x$-axis is an apartment in the maximal atlas.
However $(-1,1)\times 0\subseteq T_1$ is not even contained in an apartment of $T_1$.
\end{Num}
Nevertheless, we can do better than \ref{TopologicalRigidity}.
The completeness assumption can be removed if the dimension is high enough.
We need some more facts about the local structure of euclidean buildings
which can be found in \cite{KrWe}.
A wall in a euclidean building is called \emph{thick} if it is the
intersection of three apartments; see \cite[Sec.~10]{KrWe}.
\begin{Lem}
Let $f:X\rTo X'$ be a homeomorphism between euclidean buildings. Let $A\subseteq X$ be an
apartment and $U\subseteq A$ an open subset. Suppose that $f(U)$ is in some
apartment $A'$ in $X'$.
If $M\subseteq A$ is a thick wall, then there is a thick wall $M'\subseteq A'$ with
$f(U\cap M)=f(U)\cap M'$.

\proof
Let $p\in M\cap U$. The germ of $f(A)$ at $f(p)$ is an apartment and the
germ of $f(M)$ at $p$ is a thick wall in the spherical building
$\Delta_{f(p)}$. So $f(U\cap M)$
is locally a thick wall in $f(U)$, and therefore part of thick wall;
see \cite[10.2]{KrWe}.
\qed
\end{Lem}
A point in a euclidean building is \emph{thick} if every wall passing through
the point is thick. Thick points exist if an only if the spherical building at
infinity is thick; see \cite[10.5]{KrWe}.
\begin{Thm}
\label{TopologicalRigidityIfHigherRank}
Let $X,X'$ be irreducible euclidean buildings and suppose that
$f:X\rTo X'$ is a homeomorphism. Assume that $X$ contains a thick point
and that $X$ is not an $\RR$-tree (equivalently, that $\dim(X)\geq 2$).
Then $f$ maps apartments to apartments.

\proof
First, we note that $f$ maps the thick points in $X$ onto the thick
points in $X'$. By \cite[10.6,10.7,10.8]{KrWe} there are three possibilities:
(I) $X$ and $X'$ are euclidean cones over thick spherical buildings.
(II) $X$ and $X'$ are thick simplicial buildings.
(III) The thick points are dense in every apartment.
In case (I) and (II), both buildings are complete, so we are done by
\ref{TopologicalRigidity}. In the remaining case (III), the thick walls
are dense in every apartment $A\subseteq X$. Exactly by the same argument
as in \cite[p.~172]{BGS}, we see the following: If $U\subseteq A$ is an open 
neighborhood of $p\in A$ such that $f(U)$ is contained in an apartment
$A'\subseteq X'$, then $f:U\rTo A'$ is near $p$ an affine-linear map.
From the irreducibility of the Weyl group and the fact
that $f$ preserves thick walls near $p$ we conclude that $f$
is locally a homothety. The homothety factor is locally constant on $A$
and therefore constant on $A$. Thus $f(A)$ is isometric to $\mathbb E^n$.
By \cite[2.25]{Par}, $f(A)$ is an apartment in the maximal atlas of $X'$.
\qed
\end{Thm}
Finally, we obtain Theorem~C from the introduction. For complete
euclidean buildings with infinitely many thick points, this is
proved in \cite[Sec.~6]{KL}. We call a  euclidean building \emph{metrically
irreducible} if it is not a product of two metric spaces (with the
euclidean product of the metrics). Equivalently, $X$ is irreducible
as a building and contains a thick point, or $X=\RR$.
\begin{Thm}
\label{TopologicalRigidityInGeneral}
Let $X$ and $X'$ be metrically irreducible euclidean buildings. Assume that
$\dim(X)\geq 2$. Suppose that
$f:X\rTo X'$ is a homeomorphism. Possibly after rescaling the metric on $X'$,
there exists a unique isometry $\tilde f:X\rTo X'$ such that $\tilde f(A)=f(A)$
holds for all apartments $A\subseteq X$ of the maximal atlas of $X$.
If $X$ has more than one thick point, then $f$ has finite distance from
$\tilde f$.

\proof
We have already seen that $f$ induces a bijection between the apartments
in the maximal atlases. We use again the trichotomy from \cite[10.6,10.7,10.8]{KrWe}.
If $X$ if of type (I) or (II), this determines
the combinatorial structure of the euclidean buildings. In case (II), we see also
that $f$ has finite distance from a combinatorial isomorphism $\tilde f$.
If $X$ is of type (III),
then we proved in \ref{TopologicalRigidityIfHigherRank} that $f$ is apartmentwise
a homothety. It follows readily that $f$ is a homothety, i.e. that the metric
on $X'$ can be rescaled in such a way that $f=\tilde f$ is an isometry.
Every isometry is determined completely by its action on the apartments, so
$\tilde f$ is unique.
\qed
\end{Thm}
Using similar arguments as in \cite[6.4.3,6.4.5]{KL}, it is not difficult to extend
\ref{TopologicalRigidityInGeneral} to products of euclidean buildings and
euclidean spaces, as long as no tree factors occur. We leave this to the reader.
\begin{Num}\textbf{Remark}
Theorem~\ref{TopologicalRigidityInGeneral} is an important ingredient in proofs of quasi-isometric rigidity.
A quasi-isometry
$f:X\rTo Y$ between metric spaces always induces a homeomorphism between their asymptotic cones~%
\footnote{Asymptotic cones are truncated ultrapowers; see \cite{KrTe} for a model theoretic
viewpoint.}.
If the spaces $X$ and $Y$ are euclidean buildings or Riemannian symmetric spaces
of noncompact type, then their asymptotic cones are homeomorphic
complete euclidean buildings.
From \ref{TopologicalRigidityInGeneral} and the $\omega_1$-saturatedness of ultraproducts one 
can conclude that the quasi-isometry maps apartments (or, in the Riemannian symmetric
case, maximal flats)
Hausdorff-close to apartments; see \cite[7.1.5]{KL}. A more general result about
quasi-isometries between euclidean buildings is proved in \cite{KrWe}.
\end{Num}

\section{Euclidean buildings have finite dimension}
We recall some notions from dimension theory. See \cite{Engelking} or
\cite{Nagata} for a thorough
introduction. An \textit{open covering} $\mathfrak U$
of a space $X$ is a collection of open subsets, with $\bigcup\mathfrak U=X$.
The covering has \textit{order} $\leq n+1$ if every point $x\in X$ is contained in
at most $n+1$ elements of $\mathfrak U$. Equivalently, the dimension of the nerve
of $\mathfrak U$ is not greater than $n$.
The open covering $\mathfrak U'$ \textit{refines} $\mathfrak U$ if every 
$U'\in\mathfrak U'$ is contained in some $U\in\mathfrak U$.
The space $X$ has \textit{covering dimension} $\dim(X)\leq n$ if the following holds:
every (finite) open covering $\mathfrak U$ of $X$ has a refinement of order
$\leq n+1$. In this section we prove Theorem~B in the introduction.
For the special case $\RR$-trees, this is also proved in \cite{Burillo} and \cite{MO}
\footnote{When I wrote this paper,
I was not aware that Lang and Schlichenmaier proved in \cite[Prop.~3.3]{LaS}
that the Nagata dimension of a euclidean building with $n$-dimensional
apartments is $n$. Their result implies Theorem~\ref{DimensionThm}. Their proof
is virtually the same as ours, which follows \cite{Burillo}.}.
\begin{Thm}
\label{DimensionThm}
Let $X$ be a euclidean building with $n$-dimensional apartments. Then the covering
dimension of $X$ is
\[
\dim(X)=n.
\]
\end{Thm}
There are various other topological dimension functions, e.g. the
(small and large) \textit{inductive dimension} $\ind$ and $\Ind$.
Kleiner defined in \cite{Kl} a \textit{geometric dimension} for CAT$(\kappa)$ spaces
which we denote by $\mathrm{gdim}$. We will show that
\[
\dim(X)=\Ind(X)=\ind(X)=\mathrm{gdim}(X)=n.
\]
For metric spaces, one has always $\dim(X)=\Ind(X)\geq\ind(X)$.
If $B\subseteq X$ is a subspace, then $\ind(X)\geq\ind(B)$. 
For these facts, see \cite{Engelking,Nagata}.
Kleiner proves in \cite{Kl} that 
$\mathrm{gdim}(X)=\sup\{\dim(C)\mid C\subseteq X\text{ is compact.}\}$.
Thus $\dim(X)\geq\mathrm{gdim}(X)\geq n$.
Since $X$ contains closed sets homeomorphic to $\mathbb R^n$ and
since $\ind(\mathbb R)=\dim(\mathbb R^n)=n$, we have
$\dim(X)\geq n$.
In order to prove the theorem, we have to show that $\dim(X)\leq n$.
We use a theorem due to Burillo; see \ref{DimensionMap}.
In order to apply his result, we need that from certain
directions, euclidean buildings look 'tree-like'.
We fix an apartment $A\subseteq X$ and a regular point $\xi\in\partial A$
at infinity. Associated to these data we have the Iwasawa retraction
$\phi_{A,\xi}=\phi:X\rTo A$
which is defined as follows. If $A'\subseteq X$ is an apartment with
$\xi\in\partial A$, then $A\cap A'$ containts a Weyl chamber. Therefore
there is a unique isometry $A'\rTo A$ fixing $A\cap A'$ pointwise.
Since $X$ is the union of all apartments containing $\xi$ in their
respective boundaries, these isometries fit together to a well-defined 
$1$-Lipschitz retraction $p:X\rTo A$; see \cite[1.20]{Par}. 
We denote the fiber over $a\in A$ by $X_a=p^{-1}(a)$.
We define an ultrametric $\delta$ on $X_a$ as follows.
For $b,c\in X_a$ the geodesic rays
$(\xi,b]$ and $(\xi,c]$ branch in a point $e$,
$(\xi,b]\cap(\xi,c]=(\xi,e]$. We put 
\[
\delta(b,c)=d(e,b)+d(e,c)=2d(e,b).
\]
It is clear that this is an ultrametric on $X_a$, and that
\[
d(b,c)\leq\delta(b,c).
\]
\begin{Lem}
There exists a positive constant $L$, depending only on $\xi$, such
that 
\[
d(b,c)\leq\delta(b,c)\leq L\cdot d(b,c)
\]
holds for all $a\in A$ and all $b,c\in X_a$.

\proof
Suppose that $b\neq c$.
Let $(\xi,e]=(\xi,b]\cap(\xi,c]$ as above.
The segments $[e,b]$ and $[e,c]$ determine two different points
in the spherical building $\Sigma_eX$. These two points have the same
type, and this type depends only on $\xi$. Now
\begin{align*}
d(b,c)^2&\geq 2d(e,b)^2(1-\cos(\angle_e(b,c)))\\
&=\frac12\delta(b,c)^2(1-\cos(\angle_e(b,c)));
\end{align*}
see \cite[II.1.]{BH}.
There are only finitely many values that $\angle_e(b,c)$ can assume, since
the two points determined by $[be]$ and $[ce]$ in the spherical building
$\Sigma_eX$ are
of the same fixed type (depending only on $\xi$, not on $e$).
Moreover, $\pi\geq\angle_e(b,c)>0$. Thus $1-\cos(\angle_e(b,c))$ is bounded
away from $0$.
\qed.
\end{Lem}
Theorem~\ref{DimensionThm} follows now from Burillo's 
Theorem~\ref{DimensionMap} below.

\medskip\emph{Proof of Theorem~B.} 
By \ref{CATImpliesANR} and \ref{DimensionThm}, $X$ is a finite dimensional ANR.
An open subset of an ANR is again an ANR; see \cite[III.7.9]{Hu}.
Therefore, every open subset of a euclidean building is a finite
dimensional ANR. Since $X$ is contractible, $X$ is an AR; see \cite[III.7.2]{Hu}.
The homotopy properties follow from \ref{Mard} below.
\qed

\section{Appendix}

We first state a result about the homotopy types of ANRs. 
Appendix~1~\S2.2 in Marde\v si\'c-Segal \cite{Mar} contains more results
in this direction.
\begin{Thm}[Wall, Whitehead]\label{Mard}
Every ANR $X$ has the homotopy type of a CW complex. If $X$ has finite
covering dimension  $\dim(X)=n$, then $X$ has the homotopy type of an
$n$-dimensional CW complex, except possibly in low dimensions $n\leq 2$, where
the dimension of the CW complex may be $3$.

\proof
Every ANR $X$ is dominated by a simplicial complex $K$; see \cite[IV.6.1]{Hu}.
This means that there are continuous maps 
\[
X\rTo K\rTo X
\]
whose composite is homotopic to the identity. It follows that $X$ has the
homotopy type of a CW complex; see \cite[IV.3.8]{LW}.

Assume now that $X$ has covering dimension $n$. Let $m=\max\{3,n\}$.
The space $X$ is dominated by an $n$-dimensional simplicial complex $K$; see 
\cite[App.~1~\S2.2~Thm.~6]{Mar}. 
Then $K$ satisfies
Wall's Condition D$m$ in \cite[p.~62]{Wall}.
\begin{enumerate}
 \item[D$m$:] For all $k>m$, the groups $H_k(K)$ vanish, and $H^{m+1}(K;\mathcal B)=0$
for every coefficient bundle $\mathcal B$ on $K$.
\end{enumerate}
As Wall remarks on p.~57 in \emph{loc.cit.}, it follows that $X$ (or its CW approximation)
also satisfies condition D$m$ (since $K$ dominates $X$---all
higher dimensional homology and cohomology of $X$ has to vanish).
By Theorem E on p.~63 in \emph{loc.cit.},
$X$ has the homotopy type of an $m$-dimensional CW complex.
\qed
\end{Thm}
We remark that for $n=1$, the CW complex can also be chosen to be $1$-dimensional by Wall's
results. The $2$-dimensional case is related to the $D(2)$-problem.

The following result is proved for separable metric spaces in
\cite{Wo}. In the present form, it is stated in \cite{Gray},
but unfortunately without a proof. Our proof follows 
essentially the proof of Dugundji's Extension Theorem; see \cite[IX.6.1]{Dug}.
\begin{Thm}[Wojdys\l awski]\label{Woj}
Let $Z$ be a metric space and let $K$ be the complete simplicial complex
on a set $S$ (i.e. $K$ is the geometric realization of the poset of all finite
subsets of $S$). Suppose that there exists a continuous map $q:K\rTo Z$
with the following properties.
\begin{enumerate}
 \item[(1)] The $q$-image of the $0$-skeleton $S$ of $K$ is dense in $Z$.
 \item[(2)] If $(A_j)_{j\in\mathbb N}$ is a sequence of simplices in $K$ such that
  $(q(A_j\cap S))_{j\in\mathbb N}$ converges to $z\in Z$, then 
  $(q(A_j))_{j\in\mathbb N}$ also converges to $z$.
\end{enumerate}
Then $Z$ is an absolute extensor (AE) and in particular an absolute retract (AR)
for the class of metric spaces; see \cite[III.3.1]{Hu}.

\proof
Let $(X,d)$ be a metric space and $A\subseteq X$ closed. Let
$f:A\rTo Z$ be continuous. We have to show that $f$ extends to a continuous
map $F:X\rTo Z$. 
We proceed similarly as in \cite[IX.6.1]{Dug}. 
For each $z\in X\setminus A$ put $r_z=d(z,A)$ and $U_z=B_{r_z/2}(z)$.
These sets form an open covering $\mathfrak U$ of $X\setminus A$.
We also choose $a_z\in A$ such that $d(z,a_z)\leq 2r_z$.
Finally, we choose $s_z\in S$ such that $d(f(a_z),q(s_z))\leq r_z$
(here we use the assumption (1)).

For $x\in U_z\in\mathfrak U$ and $a\in A$ we have 
$r_z\leq d(z,a)\leq d(z,x)+d(x,a)\leq \frac12r_z+d(x,a)$, whence 
\[
r_z\leq 2d(x,a)\text{ for all }a\in A,\,x\in U_z.
\]
Moreover we have
$d(a,a_z)\leq d(a,x)+d(x,z)+d(z,a_z)\leq d(a,x)+\frac12r_z+2r_z$ and thus
\[
d(a,a_z)\leq 6d(a,x)\text{ for all }a\in A,\,x\in U_z.
\]
Let $\Phi$
be a partition of unity on $X\setminus A$ subordinate to $\mathfrak U$.
We choose for each $\phi\in\Phi$ a point $z_\phi\in X\setminus A$
such that the support of $\phi$ is in $U_\phi=U_{z_\phi}$.
Let $a_\phi=a_{z_\phi}$, $r_\phi=r_{z_\phi}$ and $s_\phi=s_{z_\phi}$.
We now define $F:X\rTo Z$ as follows. For $a\in A$ we put $F(a)=f(a)$.
For $x\in X\setminus A$, let 
\[
F(x)=q\left(\sum_\Phi \phi(x)s_\phi\right).
\]
Clearly, $F$ is continuous on $X\setminus A$. We claim that $F$ is also continuous at
each point $a\in A$. 
Let $a\in A$ and $x\in U_z$.
Then
\[
d(F(a),q(s_z))\leq d(f(a),f(a_z))+d(f(a_z),q(s_z))\leq d(f(a),f(a_z))+2d(x,a).
\]
Suppose that $(x_j)_{j\in\mathbb N}$ is a sequence in $X$ that converges to $a\in A$.
We wish to show that $\lim_j F(x_j)=F(a)$. Since $F|_A=f$ is continuous on $A$,
it suffices to consider the case where all $x_j$ are in $X\setminus A$.
Suppose that $x_j\in U_{z_j}\in\mathfrak U$. Then $\lim_jq(s_{U_{z_j}})=F(a)$ by the
inequalities above. From condition (2) and
$F(x_j)=q\left(\sum_\Phi \phi(x_j)s_\phi\right)$ we see that
$\lim_jF(x_j)=F(a)$. This shows that $F$ is continuous at $a$.
\qed
\end{Thm}
The next result shows that certain maps with ultrametric fibers cannot
lower the dimension. Our proof follows closely
\cite[Thm.~16]{Burillo}. The main difference is that we allow a Lipschitz
condition on the fibers. We use the following characterization of the covering
dimension. Recall that an open covering has \emph{order} $\leq n+1$ if its
nerve has dimension $\leq n$, or equivalently, if ever point is in at most
$n+1$ members of the covering. We say the \emph{mesh} is $\leq r$ if every member
of the covering has diameter less than $r$. The criterion that we use is the
following.
\begin{Prop}
Let $X$ be a metric space. Then $\dim(X)\leq n$ holds if an and only if
there exists a sequence of open coverings $(\mathfrak U_k)_{k\in\mathbb N}$ of
$X$ of order $\leq n+1$ and mesh $\leq r_k$, with $\lim_kr_k=0$, such that
$\mathfrak U_{k+1}$ refines $\mathfrak U_k$.

\proof
See \cite[V.1]{Nagata}.
\qed
\end{Prop}
\begin{Thm}[Burillo]\label{DimensionMap}
Let $X$ and $A$ be metric spaces and let $p:X\rTo A$ be a continuous map. We
denote the fiber over $a\in A$ by $X_a=p^{-1}(a)$. Suppose that $p$ has the
following two properties.
\begin{enumerate}
\item[(1)] Pairs of points lift isometrically:
For every $a\in A$ and $x\in X$, there exists $y\in X_a$ with
$d(x,y)=d(p(x),p(y))$.
\item[(2)] Fibers are bi-Lipschitz ultrametric:
There exists a constant $L\geq 0$ and, on each fiber $X_a$, an ultrametric 
$\delta=\delta_a:X_a\times X_a\rTo\RR$ such that 
$d(x,y)\leq\delta(x,y)\leq Ld(x,y)$ holds for all $x,y\in X_a$.
\end{enumerate}
If $\dim(A)$ is finite, then $\dim(X)\leq \dim(A)$.

\proof
Let $\mathfrak U$ be an open covering of $A$ of mesh $\leq r$ and order $\leq n+1$.
For each $U\in\mathfrak U$ we chose a point $a=a_U\in U$. For $x\in X_a$ we put
\[
W_{U,x}=\{y\in p^{-1}(U)\mid d(y,y')\leq r\text{ for some }y'\in X_a\text{ with }
\delta(y',x)\leq 3rL\}.
\]

{\em Claim 1. If $y\in W_{U,x}$ and if $z\in p^{-1}(U)$ is a point with
$d(z,y)\leq r$, then $z\in W_{U,x}$. In particular, $W_{U,x}$ is open.}

\proof
Let $y'$ be as in the definition of $W_{U,x}$ above. Using (1), choose $z'\in X_a$
with $d(z,z')\leq r$. Then $d(z',y')\leq 3r$, so $\delta(z',y')\leq 3rL$.
Therefore $\delta(x,z')\leq 3rL$ (because $\delta$ is an ultrametric).
\qed

\medskip
{\em Claim 2. For all $y,z\in W_{U,x}$ we have $d(y,z)\leq r(2+3L)$.}

\proof
Let $y,z\in W_{U,x}$ and let $y',z'\in X_a$ be points
with $d(y,y'),d(z,z')\leq r$ and $\delta(x,z'),\delta(x,y')\leq 3rL$.
Then $d(y,z)\leq r+3rL+r=r(2+3L)$.
\qed

\medskip
{\em Claim 3. Every $y\in p^{-1}(U)$ is in some $W_{U,x}$.}

\proof
Using (1), choose $x\in X_a$ with $d(x,y)=d(p(x),a)\leq r$. Then $y\in W_{U,x}$.
\qed

\medskip
{\em Claim 4. If $W_{U,x}\cap W_{U,z}\neq\emptyset$, then
$W_{U,x}=W_{U,z}$.}

\proof
Suppose that $y\in W_{U,x}\cap W_{U,z}$. Let $y'\in X_a$ with
$d(y,y')\leq r$. By Claim 1 we have $y'\in W_{U,x}\cap W_{U,z}$.
Since $\delta$ is ultrametric, we have $\delta(x,z)\leq 3rL$ and therefore
$W_{U,x}=W_{U,z}$.
\qed

\medskip
{\em Claim 5. The set $\mathfrak W=\{W_{U,x}\mid U\in\mathfrak U\text{ and }
p(x)=a_U\}$ is an open covering of $X$ of order $\leq n+1$ and mesh $\leq r(2+3L)$.}

\proof
Only the claim about the order remains to be shown. 
Suppose that $W_{U_0,x_0},\ldots,W_{U_m,x_m}$ are pairwise
distinct, with $W_{U_0,x_0}\cap\cdots\cap W_{U_m,x_m}\neq \emptyset$.
By Claim 4, we have $U_i\neq U_j$ for $i\neq j$. Moreover,
$U_0\cap\cdots\cap U_m\neq \emptyset$. It follows that $m\leq n$.
\qed

\medskip
{\em Claim 6. Let $\tilde{\mathfrak U}$ be an open covering of $A$
of mesh $\leq \tilde r\leq r/(2+3L)$ and order $\leq n+1$ that refines $\mathfrak U$.
Let $\tilde{\mathfrak W}$ be constructed from $\tilde{\mathfrak U}$ 
and $\tilde r$ as above. Then $\tilde{\mathfrak W}$ refines $\mathfrak W$.}

\proof
Let $\tilde W_{V,v}\in\tilde {\mathfrak W}$.
We choose $U\in\mathfrak U$ with $V\subseteq U$ and $x\in X_{a_U}$ such that
$v\in W_{U,x}$. For $y\in \tilde W_{V,v}$ we have (by Claim 2) that
$d(y,v)\leq \tilde r(2+3L)\leq r$.
By Claim 1 we have $y\in W_{U,x}$.
\qed

\medskip
The claim of the theorem now follows. We choose a sequence of coverings $\mathfrak U_k$
of $A$ of order $\leq n+1$ and mesh $\leq r_k$, with
$(2+3L)r_{k+1}\leq r_k$. The resulting sequence of coverings $\mathfrak W_k$ of $X$
then has order $\leq n+1$ and mesh $\leq r_k(2+3L)$.
\qed
\end{Thm}

\raggedright
Linus Kramer\\
Mathematisches Institut, 
Universit\"at M\"unster,
Einsteinstr. 62,
48149 M\"unster,
Germany\\
\makeatletter
e-mail: {\tt linus.kramer{@}uni-muenster.de}
\end{document}